\documentclass[11pt]{article}
\usepackage{latexsym,amsmath,amsfonts,graphicx}
\DeclareSymbolFont{lettersA}{U}{pxmia}{m}{it}
\DeclareMathSymbol{\piup}{\mathord}{lettersA}{"19}
\begin{document}
\begin{center}
{\LARGE\textbf{Focusing of Maximum Vertex Degrees\\
in Random Faulty Scaled Sector Graphs}}\\
\bigskip
Yilun Shang

Institute for Cyber Security, University of Texas at San Antonio\\
San Antonio, Texas 78249, USA

\ttfamily{shylmath@hotmail.com}

\end{center}

\begin{abstract}
In this paper we study the behavior of maximum out/in-degree of
binomial/Poisson random scaled sector graphs in the presence of
random vertex and edge faults. We prove that the probability
distribution of maximum degrees for random faulty scaled sector
graphs with $n$ vertices, where each vertex spans a sector of
$\alpha$ radians, with radius $r_n\ll \sqrt{\ln n/n}$, becomes
concentrated on two consecutive integers, under some natural
assumptions of fault probabilities.
\bigskip

\textbf{2010 Mathematics Subject Classification:} 05C80, 05C20,
60C05

\smallskip
\textbf{Keywords:} Random scaled sector graph; Random geometric
graph; Stein-Chen method; Degree; Percolation.

\end{abstract}

\bigskip
\normalsize

\section{Introduction}

Wireless ad-hoc communication networks of sensors gain increasing
importance in telecommunication society during the last decades
\cite{11}. The rapid developments in peer-to-peer networks have
resulted in the emergence of various models including random
distance graph and geometric graph models \cite{20,41,42}. Besides,
the random scaled sector graph model proposed in \cite{1} aims to
provide a tool for the analysis of routing and transmission of
information in sensor networks communicating through optical devices
or directional antennae. In this setting, a large number of randomly
scattered transmitters are located in some geographical area. A
transmitter can orient its laser beam within a fixed radius $r$ of
it and also in a range of directions. In other words, the scanning
area of a transmitter can be viewed as a sector of radius $r$ and
radian $\alpha$ (see Fig. 1 below). In practical applications,
various hostile unexpected environments are inevitable and thus
should be taken into account. Every transmitter may have a failure
probability understood as the sensor becoming inoperative due to
mechanical damages or power drain. This is usually described as site
percolation. Every connection may also have a failure probability
understood as the failure of communication because of bad weather or
terrain obstacles. This is usually described as bond percolation. If
we think of the sensor network as a graph with each transmitter as a
vertex and each (directed) connection as an arc, we will have a
faulty directed graph describing the above situation. We shall refer
to this model as random faulty scaled sector graph (see Section 2
for a formal definition), and some recent results regarding its
degree distributions, small subgraph and coloring problems can be
found in \cite{30,31,40}.

In this paper, we investigate extreme degrees of the above mentioned
random digraph model. A concentration result (Theorem 1) shows that
the maximum out-/in-degrees of random faulty scaled sector graph are
almost determined under some natural assumptions of fault
probabilities (see below for details). A similar focusing phenomenon
has been discovered in classical Erd\"os-R\'enyi random graph
theory, see e.g.\cite[Chap. 3]{12}. Theorem 1 adds to the asymptotic
bound of maximum degree of \cite{1} in the thermodynamic and
sub-connective regimes by including Poisson point process case. It
extends the maximum degree focusing result in \cite[Chap. 6]{6} by
considering digraphs and including fault probabilities. It also
offers a partial answer to the open problems suggested in \cite{2}.

Poisson approximation by Stein-Chen method is used here as in
\cite{7}, where Penrose incorporated the geometric clique number
with scan statistic via a ``clustering rule'' giving a concentration
result. The idea behind can be traced back to \cite{9}, for
instance. We mention that if the clustering rule $h$ is properly
chosen, the geometric maximum degree is also contained in that
framework. In fact, let $h(\mathcal{X})=0$ if $\mathcal{X}$ is not
contained in some ball $B(x,r)$ and otherwise let $h(\mathcal{X})$
be equal to the maximum degree$+1$ of the geometric subgraph induced
by $\mathcal{X}$. It is straightforward to verify that $h$ satisfies
the requirements in \cite{7}.

Compared with the more recent random scaled sector graph model,
random geometric graph models are widely-studied, see e.g.
\cite{38,37,36,35,34,6,33,32,40} and references therein, in which
connections are isotropic and thus undirected. Some properties such
as connectivity and layout problems of geometric graphs featuring
fault probability have been addressed, see e.g. \cite{2,16,3,8}.
Algorithmic aspects are also dealt with in various contexts, see
e.g. \cite{39,4,13}. From a percolation analysis point of view, some
relevant results coherent with ours can be found in \cite{5}.

\section{Statement of main results}

\begin{figure}[hbt]
\centering
\scalebox{0.6}{\includegraphics[245pt,117pt][520pt,352pt]{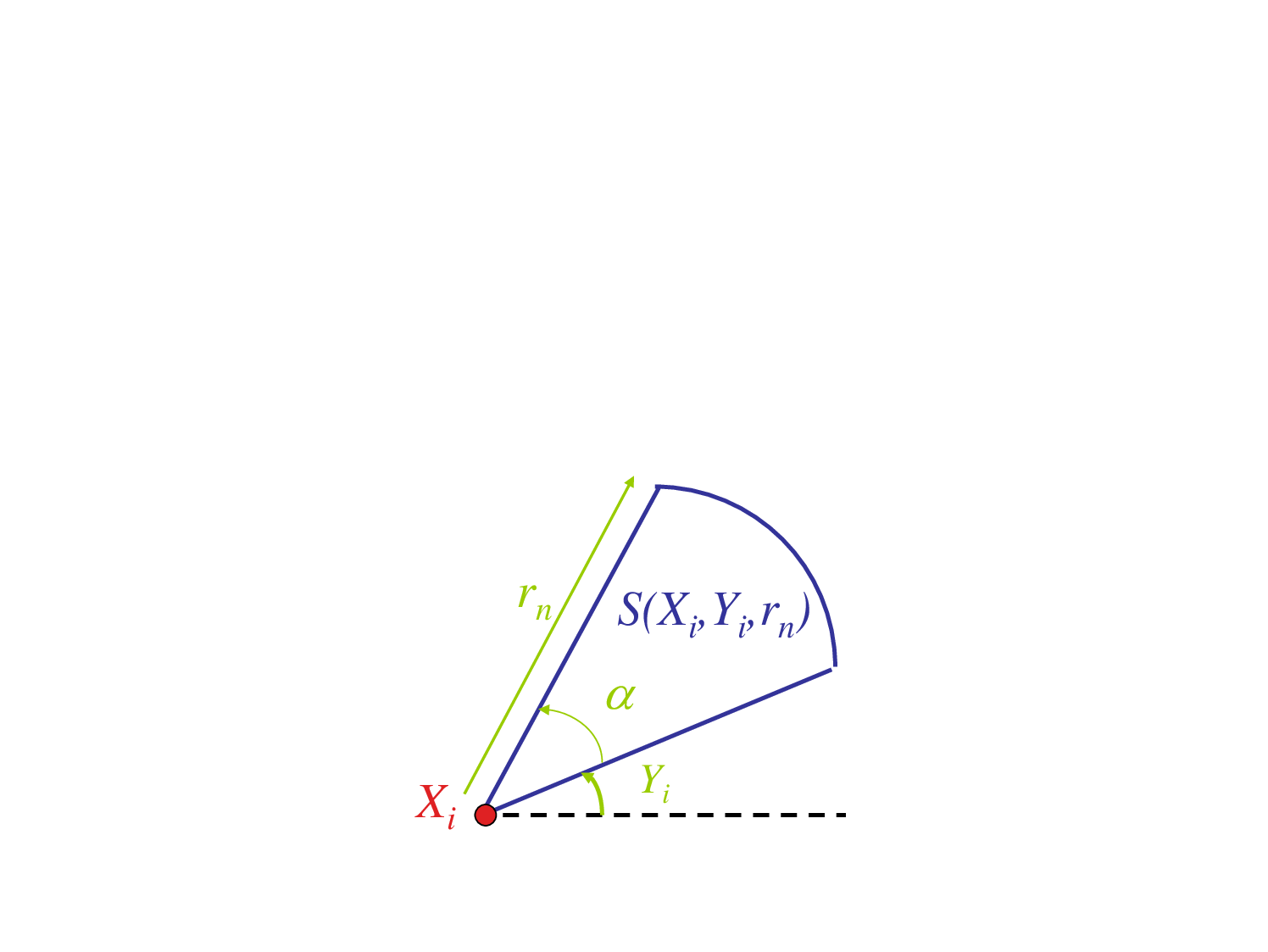}}
\caption{An illustration of a sector $S(X_i,Y_i,r_n)$.}
\end{figure}

To begin with, we fix some notations that will be used in the
derivation of the main results. Given a sequence
$\mathcal{X}_n=\{X_1,X_2,\cdots,X_n\}$ of independently and
uniformly distributed $(i.u.d.)$ random points in the square area
$[0,1]^2$ with common density function $f=1_{[0,1]^2}$. Here, for a
set $A\subseteq\mathbb{R}^2$, $1_A$ is its indicator. We equip
$\mathbb{R}^2$ with Euclidean norm and let $\theta$ be the area of
unit disk. Notice that $\theta=\piup$ in this context, and we
therefore use only $\piup$ in the sequel to simplify notations. Let
$\mathcal{Y}_n=\{Y_1,Y_2,\cdots,Y_n\}$ be $i.u.d.$ random variables
taking values in $[0,2\piup)$ and $\alpha\in(0,2\piup]$ be fixed.
Associate every point $X_i\in \mathcal{X}_n$ a sector, which is
centered at $X_i$, with radius $r_n$, central angle $\alpha$ and
elevation $Y_i$ with respect to the horizontal direction
anticlockwise. This sector is denoted as $S(X_i,Y_i,r_n)$; see Fig.
1. We denote by $G_\alpha(\mathcal{X}_n,\mathcal{Y}_n,r_n)$ the
digraph with vertex set $\mathcal{X}_n$, and with an arc
$(X_i,X_j)$, $i\not=j$, presents if and only if $X_j\in
S(X_i,Y_i,r_n)$. For any $A\subseteq\mathbb{R}^2$, let
$\mathcal{X}_n(A)$ denote the number of vertices in
$\mathcal{X}_n\cap A$. Given $\lambda>0$, denote by $Poi(\lambda)$
the Poisson distribution with parameter $\lambda$. A random variable
$X\sim Poi(\lambda)$ represents that $X$ obeys the Poisson
distribution. The usual Poisson version
$G(\mathcal{P}_n,\mathcal{Y}_{N_n},r_n)$ is defined similarly, where
$\mathcal{P}_n=\{X_1,X_2,\cdots,X_{N_n}\}$ and $N_n\sim Poi(n)$.

For every vertex $X_i$, we assign a failure probability $v_n$,
independent of other nodes and the point process. If a vertex is
failed, the vertex itself together with all out/in-edges adjacent to
it is then removed from the graph. Given the presence of vertex
$X_i$, we associate every out-edge with a failure probability $q_n$
independently. This random faulty scaled sector graph is thus
denoted by $G_{\alpha}(\mathcal{X}_n,\mathcal{Y}_n,v_n,q_n,r_n)$.
Likewise, we can define the Poisson version
$G_{\alpha}(\mathcal{P}_n,\mathcal{Y}_{N_n},v_n,q_n,r_n)$. Let
$\Delta_n^{out}$, $\Delta_n^{in}$ be the maximum out-/in-degree of
$G_{\alpha}(\mathcal{X}_n,\mathcal{Y}_n,v_n,q_n,r_n)$ respectively,
and $\Delta_n^{'out}$, $\Delta_n^{'in}$ be the maximum
out-/in-degree of
$G_{\alpha}(\mathcal{P}_n,\mathcal{Y}_{N_n},v_n,q_n,r_n)$
respectively. Our main result is the following:

\medskip
\noindent\textbf{Theorem 1.}\itshape \quad Suppose $v_n\rightarrow
v\in[0,1)$ and $q_n\rightarrow q\in[0,1)$, as $n\rightarrow\infty$.
Suppose $\mu_n:=\frac{\alpha}{2}nr_n^2(1-v_n)(1-q_n)$, and that
$\inf_{n>0}\mu_n>0$, and that $\mu_n^{1+\varepsilon}=o(\ln n)$ for
some $\varepsilon>0$. Then there exists a sequence
$\{k_n\}_{n\ge1}$, set $\xi_n=P(Poi(\mu_n)\ge k_n)$, such that we
have\\
$$P(\Delta_n^{'out}=k_n-1)-e^{-n(1-v_n)\xi_n}\rightarrow0,$$
$$P(\Delta_n^{'out}=k_n)+e^{-n(1-v_n)\xi_n}\rightarrow1,$$ as
$n\rightarrow\infty$. The same thing holds for $\Delta_n^{'in}$,
$\Delta_n^{out}$ and $\Delta_n^{in}$, respectively.\normalfont

\section{Proofs}

To prove the asymptotic focusing phenomenon, we first give a general
non-asymptotic Poisson approximation lemma, which may be useful in
some other cases. Let $W_{j,n}^{out}(r)$, $W_{j,n}^{in}(r)$ be the
number of vertices of out-/in-degree $j$ in
$G_{\alpha}(\mathcal{X}_n,\mathcal{Y}_n,v,q,r)$, respectively. Let
$W_{j,\lambda}^{'out}(r)$, $W_{j,\lambda}^{'in}(r)$ be the number of
vertices of out-/in-degree $j$ in
$G_{\alpha}(\mathcal{P}_{\lambda},\mathcal{Y}_{N_{\lambda}},v,q,r)$,
respectively. For $A\subseteq\mathbb{N}\cup\{0\}$, set
$W_{A,\lambda}^{'out}(r):=\sum_{j\in A}W_{j,\lambda}^{'out}(r)$ and
$W_{A,\lambda}^{'in}(r):=\sum_{j\in A}W_{j,\lambda}^{'in}(r)$. The
total variation distance between the laws of non-negative integer
valued random variables $X$, $Y$ is defined by
$$
d_{TV}(X,Y)=\sup_{A\subseteq\mathbb{N}\cup\{0\}}\{|P(X\in A)-P(Y\in A)|\}.
$$
Given $x\in\mathbb{R}^2$ and $r>0$, define $B(x,r)$ the disk with
center $x$ and radius $r$. Let $c$, $c'$ be various positive
constants throughout the paper, and the values may change from line
to line.

\medskip
\noindent\textbf{Lemma 1.}\itshape \quad Suppose a density function
$g$ is continuous a.e. on $\mathbb{R}^2$. Let $0\le v, q<1$, $r,
\lambda>0$ and $A\subseteq\mathbb{N}\cup\{0\}$. Then,
$$d_{TV}(W_{A,\lambda}^{'out}, Poi(EW_{A,\lambda}^{'out}))\le \min\Big(1,\frac{1}{EW_{A,\lambda}^{'out}}\Big)\cdot(I_1^{out}+I_2^{out})$$
where,
\begin{multline}
I_1^{out}:=\frac{(1-v)^2\lambda^2}{4\piup^2}\int_{\mathbb{R}^2}\int_0^{2\piup}P[\mathcal{P}_{\lambda(1-q)(1-v)}(S(x_1,y_1,r))\in
A]\mathrm{d}y_1\\
\int_{B(x_1,3r)}\int_0^{2\piup}P[\mathcal{P}_{\lambda(1-q)(1-v)}(S(x_2,y_2,r))\in
A]\mathrm{d}y_2g(x_2)\mathrm{d}x_2g(x_1)\mathrm{d}x_1,\nonumber
\end{multline}
and
\begin{multline}
I_2^{out}:=\frac{(1-v)^2\lambda^2}{4\piup^2}\\
\cdot\int_{\mathbb{R}^2}\int_{B(x_1,3r)}\int_0^{2\piup}\int_0^{2\piup}P\big[\{\mathcal{P}_{\lambda(1-q)(1-v)}(S(x_1,y_1,r))
+1_{[x_2\in S(x_1,y_1,r)]\cap[(x_1,x_2)\ not\ fails]}\in A\}\\
\cap
\{\mathcal{P}_{\lambda(1-q)(1-v)}(S(x_2,y_2,r))+1_{[x_1\in
S(x_2,y_2,r)]\cap[(x_2,x_1)\ not\ fails]}\in
A\}\big]\\
\mathrm{d}y_1\mathrm{d}y_2g(x_2)\mathrm{d}x_2g(x_1)\mathrm{d}x_1.\nonumber
\end{multline}
Likewise,
$$d_{TV}(W_{A,\lambda}^{'in}, Poi(EW_{A,\lambda}^{'in}))\le \min\Big(1,\frac{1}{EW_{A,\lambda}^{'in}}\Big)(I_1^{in}+I_2^{in})$$
where,
\begin{multline}
I_1^{in}:=(1-v)^2\lambda^2\int_{\mathbb{R}^2}P[\mathcal{\tilde{P}}_{\lambda(1-q)(1-v)}(B(x_1,r))\in A]\\
\int_{B(x_1,3r)}P[\mathcal{\tilde{P}}_{\lambda(1-q)(1-v)}(B(x_2,r))\in
A]g(x_2)\mathrm{d}x_2g(x_1)\mathrm{d}x_1,\nonumber
\end{multline}
and
\begin{multline}
I_2^{in}:=(1-v)^2\lambda^2\\
\cdot\int_{\mathbb{R}^2}\int_{B(x_1,3r)}P\big[\{\mathcal{\tilde{P}}_{\lambda(1-q)(1-v)}(B(x_1,r))+1_{[x_1\in
S(x_2,y_2,r)]\cap[(x_2,x_1)\ not\ fails]}\in
A\}\\
\cap\{\mathcal{\tilde{P}}_{\lambda(1-q)(1-v)}(B(x_2,r))+1_{[x_2\in
S(x_1,y_1,r)]\cap[(x_1,x_2)\ not\ fails]}\in
A\}\big]g(x_2)\mathrm{d}x_2g(x_1)\mathrm{d}x_1. \nonumber
\end{multline}
Here, we denote by $\mathcal{\tilde{P}}_{\lambda}$ a Poisson point
process with intensity $(\lambda\alpha/2\piup)g$, which is the
thinning of $\mathcal{P}_{\lambda}$ whose intensity is $\lambda g$.
\normalfont

\medskip
\noindent\textbf{Proof}. Given $m\in\mathbb{N}$, partition
$\mathbb{R}^2$ into squares of side $2^{-m}$, with the origin lies
at a square corner. Label these squares as $D_{m,1},D_{m,2},\cdots$,
and denote the center of $D_{m,i}$ as $a_{m,i}$. For each
$x\in\mathbb{R}^2$ and for each $m$, $i$, define $Y_x$, $Y_{m,i}$ as
independent copies of $Y_1$.

For out-degree, set
$$\xi_{m,i}:=1_{[\mathcal{P}_{\lambda(1-v)}(D_{m,i})=1]\cap[\mathcal{P}_{\lambda(1-q)(1-v)}(S(a_{m,i},Y_{m,i},r)\backslash D_{m,i})\in A]}$$
Set $p_{m,i}:=E\xi_{m,i}$, $p_{m,i,j}:=E[\xi_{m,i}\xi_{m,j}]$.
Define an adjacency relation $\sim_m$ on $\mathbb{N}$ by putting
$i\sim_m j$ if and only if $0<||a_{m,i}-a_{m,j}||\le 3r$, and define
the corresponding adjacency neighborhood
$\mathcal{N}_{m,i}:=\{j\in\mathbb{N}|\ ||a_{m,i}-a_{m,j}||\le 3r\}$.
Let $Q_n:=[-n,n]^2$ and
$\mathcal{I}_{m,n}:=\{i\in\mathbb{N}|D_{m,i}\subseteq Q_n\}$. Set
$\mathcal{N}_{m,n,i}:=\mathcal{N}_{m,i}\cap\mathcal{I}_{m,n}$. Thus
$(\mathcal{I}_{m,n},\sim_m)$ is a dependency graph for random
variables $\xi_{m,i}$, $i\in\mathcal{I}_{m,n}$.

Define
$\tilde{W}_{m,n}^{out}:=\sum_{i\in\mathcal{I}_{m,n}}\xi_{m,i}$, then
we observe that
$W_{A,\lambda}^{'out}=\lim_{n\rightarrow\infty}\lim_{m\rightarrow\infty}\tilde{W}_{m,n}^{out}$.
By Theorem 1 of \cite{14},
\begin{equation}
d_{TV}(\tilde{W}_{m,n}^{out},Poi(E\tilde{W}_{m,n}^{out}))\le
\min\Big(1,\frac{1}{E\tilde{W}_{m,n}^{out}}\Big)(a_1(m,n)+a_2(m,n))\label{1}
\end{equation}
where
$$a_1(m,n):=\sum_{i\in\mathcal{I}_{m,n}}\sum_{j\in\mathcal{N}_{m,n,i}}p_{m,i}p_{m,j},\qquad a_2(m,n):=\sum_{i\in\mathcal{I}_{m,n}}\sum_{j\in\mathcal{N}_{m,n,i}\backslash\{i\}}p_{m,i,j}.$$
Define $w_m(x):=2^{2m}p_{m,i}1_{[x\in D_{m,i}]}$, wherefore
$\int_{Q_n}w_m(x)\mathrm{d}x=\sum_{i\in\mathcal{I}_{m,n}}p_{m,i}$.
If $f$ is continuous at $x$, we have
$\lim_{m\rightarrow\infty}w_m(x)=(1-v)\lambda
g(x)P[\mathcal{P}_{\lambda(1-q)(1-v)}(S(x,Y_x,r))\in A]$, by the
mean value theorem. Observe that $w_m(x)\le
2^{2m}E\mathcal{P}_{\lambda(1-v)}(D_{m,i})\le \lambda g_{\max}$, so
by the dominated convergence theorem, we obtain
$$\lim_{m\rightarrow\infty}E\tilde{W}_{m,n}^{out}=(1-v)\lambda\int_{Q_n}P[\mathcal{P}_{\lambda(1-q)(1-v)}(S(x,Y_x,r))\in A]g(x)\mathrm{d}x$$
and by the Fubini theorem and Palm theory for (marked) Poisson point
process \cite{6}(Sect. 1.7), we have
\begin{eqnarray*}
\lim_{n\rightarrow\infty}\lim_{m\rightarrow\infty}E\tilde{W}_{m,n}^{out}&\hspace{-7pt}=&\hspace{-7pt}(1-v)\lambda\int_{\mathbb{R}^2}P[\mathcal{P}_{\lambda(1-q)(1-v)}(S(x,Y_x,r))\in
A]g(x)\mathrm{d}x\\
&\hspace{-7pt}=&\hspace{-7pt}\frac{(1-v)\lambda}{2\piup}\int_{\mathbb{R}^2}\int_0^{2\piup}P[Poi\big(\int_{S(x,y,r)}\lambda(1-v)(1-q)g(z)\mathrm{d}z\big)\in
A]g(x)\mathrm{d}y\mathrm{d}x\\
&\hspace{-7pt}=&\hspace{-7pt}EW_{A,\lambda}^{'out}
\end{eqnarray*}
For $x_1\in D_{m,i}$, $x_2\in D_{m,j}$, define
$u_m(x_1,x_2):=2^{4m}p_{m,i}p_{m,j}1_{[j\in \mathcal{N}_{m,i}]}$ and
$v_m(x_1,x_2):=2^{4m}p_{m,i,j}1_{[j\in
\mathcal{N}_{m,i}\backslash\{i\}]}$. Therefore, we have
$a_1(m,n)=\int_{Q_n}\int_{Q_n}u_m(x_1,x_2)\mathrm{d}x_1\mathrm{d}x_2$
and
$a_2(m,n)\\=\int_{Q_n}\int_{Q_n}v_m(x_1,x_2)\mathrm{d}x_1\mathrm{d}x_2$.
For different continuous points $x_1$, $x_2$ of $g$, if also
$||x_1-x_2||\not=r$ and $||x_1-x_2||\not=3r$, then
\begin{multline}
\lim_{m\rightarrow\infty}u_m(x_1,x_2)=\frac{(1-v)^2\lambda^2}{4\piup^2}g(x_1)g(x_2)\int_0^{2\piup}P[\mathcal{P}_{\lambda(1-q)(1-v)}(S(x_1,y_1,r))\in
A]\mathrm{d}y_1\\
\cdot\int_0^{2\piup}P[\mathcal{P}_{\lambda(1-q)(1-v)}(S(x_2,y_2,r))\in
A]\mathrm{d}y_2\cdot1_{[B(x_1,3r)]}(x_2)\nonumber
\end{multline}
Similarly,
\begin{multline}
\lim_{m\rightarrow\infty}v_m(x_1,x_2)=\frac{(1-v)^2\lambda^2}{4\piup^2}g(x_1)g(x_2)\\
\cdot\int_0^{2\piup}\int_0^{2\piup}P\big[\{\mathcal{P}_{\lambda(1-q)(1-v)}(S(x_1,y_1,r))
+1_{[x_2\in S(x_1,y_1,r)]\cap[(x_1,x_2)\ not\ fails]}\in A\}\\\cap
\{\mathcal{P}_{\lambda(1-q)(1-v)}(S(x_2,y_2,r))+1_{[x_1\in
S(x_2,y_2,r)]\cap[(x_2,x_1)\ not\ fails]}\in
A\}\big]\mathrm{d}y_1\mathrm{d}y_2\cdot1_{[B(x_1,3r)]}(x_2)\nonumber
\end{multline}
For $x_1\in D_{m,i}$, $x_2\in D_{m,j}$, we have
$$u_m(x_1,x_2)\le
2^{4m}E\mathcal{P}_{\lambda(1-v)}(D_{m,i})E\mathcal{P}_{\lambda(1-v)}(D_{m,j})\le
g_{\max}^2$$ and $$v_m(x_1,x_2)\le
2^{4m}E\mathcal{P}_{\lambda(1-v)}(D_{m,i})E\mathcal{P}_{\lambda(1-v)}(D_{m,j})1_{[i\not=j]}\le
g_{\max}^2.$$ Hence, by the dominated convergence theorem we have
$\lim_{n\rightarrow\infty}\lim_{m\rightarrow\infty}a_1(m,n)=I_1^{out}$
and
$\lim_{n\rightarrow\infty}\lim_{m\rightarrow\infty}a_2(m,n)=I_2^{out}$.
Hence, the out-degree case is proved by taking limit in both sides
of (\ref{1}).

For in-degree, set
$$\eta_{m,i}:=1_{[\mathcal{P}_{\lambda(1-v)}(D_{m,i})=1]\cap[\#\{x'\in\mathcal{P}_{\lambda(1-q)(1-v)}|D_{m,i}\subseteq S(x',Y_{x'},r)\}\in A]}$$
Set $q_{m,i}:=E\eta_{m,i}$, $q_{m,i,j}:=E[\eta_{m,i}\eta_{m,j}]$.
Define the dependency graph for random variables $\eta_{m,i}$,
$i\in\mathcal{I}_{m,n}$ just as above.

Define
$\tilde{W}_{m,n}^{in}:=\sum_{i\in\mathcal{I}_{m,n}}\eta_{m,i}$, then
we observe that
$W_{A,\lambda}^{'in}=\lim_{n\rightarrow\infty}\lim_{m\rightarrow\infty}\tilde{W}_{m,n}^{in}$.
By Theorem 1 of \cite{14},
\begin{equation}
d_{TV}(\tilde{W}_{m,n}^{in},Poi(E\tilde{W}_{m,n}^{in}))\le
\min\Big(1,\frac{1}{E\tilde{W}_{m,n}^{in}}\Big)(b_1(m,n)+b_2(m,n))\label{2}
\end{equation}
where
$$b_1(m,n):=\sum_{i\in\mathcal{I}_{m,n}}\sum_{j\in\mathcal{N}_{m,n,i}}q_{m,i}q_{m,j},\qquad b_2(m,n):=\sum_{i\in\mathcal{I}_{m,n}}\sum_{j\in\mathcal{N}_{m,n,i}\backslash\{i\}}q_{m,i,j}.$$
Reset $w_m(x):=2^{2m}q_{m,i}1_{[x\in D_{m,i}]}$, then
$\int_{Q_n}w_m(x)\mathrm{d}x=\sum_{i\in\mathcal{I}_{m,n}}q_{m,i}$.
If $f$ is continuous at $x$, we have
$\lim_{m\rightarrow\infty}w_m(x)=(1-v)\lambda
g(x)P[\#\{x'\in\mathcal{P}_{\lambda(1-q)(1-v)}|x\in
S(x',Y_{x'},r)\}\in A]$, by mean value theorem for integrals.
Observe that $w_m(x)\le
2^{2m}E\mathcal{P}_{\lambda(1-v)}(D_{m,i})\le \lambda g_{\max}$, so
by the dominated convergence theorem, we obtain
$$\lim_{m\rightarrow\infty}E\tilde{W}_{m,n}^{in}=(1-v)\lambda\int_{Q_n}P[\#\{x'\in\mathcal{P}_{\lambda(1-q)(1-v)}|x\in
S(x',Y_{x'},r)\}\in A]g(x)\mathrm{d}x$$ and
\begin{eqnarray*}
\lim_{n\rightarrow\infty}\lim_{m\rightarrow\infty}E\tilde{W}_{m,n}^{in}&\hspace{-7pt}=&\hspace{-7pt}(1-v)\lambda\int_{\mathbb{R}^2}P[Poi\big(\frac{\lambda\alpha(1-q)(1-v)}{2\piup}\int_{B(x,r)}g(z)\mathrm{d}z\big)\in A]g(x)\mathrm{d}x\\
&\hspace{-7pt}=&\hspace{-7pt}EW_{A,\lambda}^{'in}
\end{eqnarray*}
For $x_1\in D_{m,i}$, $x_2\in D_{m,j}$, reset
$u_m(x_1,x_2):=2^{4m}q_{m,i}q_{m,j}1_{[j\in \mathcal{N}_{m,i}]}$ and
$v_m(x_1,x_2):=2^{4m}q_{m,i,j}1_{[j\in
\mathcal{N}_{m,i}\backslash\{i\}]}$. Therefore, we have
$b_1(m,n)=\int_{Q_n}\int_{Q_n}u_m(x_1,x_2)\mathrm{d}x_1\mathrm{d}x_2$
and
$b_2(m,n)\\=\int_{Q_n}\int_{Q_n}v_m(x_1,x_2)\mathrm{d}x_1\mathrm{d}x_2$.
For different continuous points $x_1$, $x_2$ of $g$, if also
$||x_1-x_2||\not=r$ and $||x_1-x_2||\not=3r$, then
\begin{multline}
\lim_{m\rightarrow\infty}u_m(x_1,x_2)=(1-v)^2\lambda^2
g(x_1)g(x_2)P[\mathcal{\tilde{P}}_{\lambda(1-q)(1-v)}(B(x_1,r))\in
A]\\
\cdot P[\mathcal{\tilde{P}}_{\lambda(1-q)(1-v)}(B(x_2,r))\in
A]1_{[B(x_1,3r)]}(x_2)\nonumber
\end{multline}
Similarly,
\begin{multline}
\lim_{m\rightarrow\infty}v_m(x_1,x_2)=(1-v)^2\lambda^2g(x_1)g(x_2)\\
\cdot P[\{\mathcal{\tilde{P}}_{\lambda(1-q)(1-v)}(B(x_1,r))
+1_{[x_1\in S(x_2,Y_2,r)]\cap[(x_2,x_1)\ not\ fails]}\in A\}\\\cap
\{\mathcal{\tilde{P}}_{\lambda(1-q)(1-v)}(B(x_2,r))+1_{[x_2\in
S(x_1,Y_1,r)]\cap[(x_1,x_2)\ not\ fails]}\in A\}]\cdot
1_{[B(x_1,3r)]}(x_2)\nonumber
\end{multline}
By similar arguments in out-degree case, we have
$\lim_{n\rightarrow\infty}\lim_{m\rightarrow\infty}b_1(m,n)=I_1^{in}$
and
$\lim_{n\rightarrow\infty}\lim_{m\rightarrow\infty}b_2(m,n)=I_2^{in}$.
Hence, we conclude the proof by taking limit in both sides of
(\ref{2}). $\Box$
\medskip

For cleanness of the expressions, we will shift our battlefield from
$[0,1]^2$ to $[-1/2,1/2]^2$. From now on, we take $f=1_Q$,
$Q:=[-1/2,1/2]^2$. Let $\mathcal{H}_{\lambda}$ be the homogeneous
Poisson point process with intensity $\lambda$ on $\mathbb{R}^2$ and
$|\cdot|$ be Lebesgue measure.

\medskip
\noindent\textbf{Proposition 1.}\itshape \quad Let
$\mu_n:=\frac{\alpha}{2}nr_n^2(1-v_n)(1-q_n)$, and suppose
$\lim_{n\rightarrow\infty}r_n=0$ and $\inf_{n>0}\mu_n>0$. Suppose
$\{j_n\}_{n\ge 1}$ is an $\mathbb{N}-$valued sequence such that for
some $\varepsilon>0$,
\begin{equation}
\lim_{n\rightarrow\infty}j_n/\mu_n^{1+\varepsilon}=\infty \label{3}
\end{equation}
Set $\xi_n:=P(Poi(\mu_n)\ge j_n)$. Then
$$\lim_{n\rightarrow\infty}[P(\Delta_n^{'out}<j_n)-e^{-n(1-v_n)\xi_n}]=0$$
and
$$\lim_{n\rightarrow\infty}[P(\Delta_n^{'in}<j_n)-e^{-n(1-v_n)\xi_n}]=0.$$\normalfont

\medskip
\noindent\textbf{Proof}. For out-degree, set
$$W_n^{'out}:=\sum_{i=1}^{N_{n(1-v_n)}}1_{[\mathcal{P}_{n(1-v_n)(1-q_n)}(S(X_i,Y_i,r_n))\ge j_n+1]\cap[X_i\in
Q]}$$ Then by Palm theory for (marked) Poisson process,
$$EW_n^{'out}\sim (1-v_n)n\int_QP[Poi\big(\frac{\alpha}{2}nr_n^2(1-v_n)(1-q_n)\big)\ge j_n]\mathrm{d}x=(1-v_n)n\xi_n,$$
as $n\rightarrow\infty$. Now take $A=\mathbb{Z}\cap[j_n,\infty)$,
$\lambda=n$, $r=r_n$, $v=v_n$, $q=q_n$, $g=f$ in Lemma 1, we then
obtain
\begin{equation}
\big|P(W_n^{'out}=0)-e^{-EW_n^{'out}}\big|\le
\min\Big(1,\frac{2}{(1-v_n)n\xi_n}\Big)(I_{1,n}^{out}+I_{2,n}^{out}),\label{4}
\end{equation}
for large enough $n$. We have
\begin{eqnarray*}
I_{1,n}^{out}&=&\frac{(1-v_n)^2n^2}{4\piup^2}\int_{\mathbb{R}^2}\int_0^{2\piup}P[\mathcal{P}_{n(1-q_n)(1-v_n)}(S(x_1,y_1,r_n))\ge
j_n]\mathrm{d}y_1\\
&
&\int_{B(x_1,3r_n)}\int_0^{2\piup}P[\mathcal{P}_{n(1-q_n)(1-v_n)}(S(x_2,y_2,r_n))\ge
j_n]\mathrm{d}y_2f(x_2)\mathrm{d}x_2f(x_1)\mathrm{d}x_1\\
&\le&n^2\xi_n^2\piup(3r_n)^2
\end{eqnarray*}
Therefore, by (\ref{3}), $((1-v_n)n\xi_n)^{-1}I_{1,n}^{out}\le
c\mu_n\xi_n\rightarrow0.$ On the other hand, we have,
\begin{eqnarray*}
I_{2,n}^{out}&=&\frac{(1-v_n)^2n^2}{4\piup^2}\int_Q\int_{Q\cap
B(x_1,3r_n)}\int_0^{2\piup}\int_0^{2\piup}P\big[\{\mathcal{P}_{n(1-q_n)(1-v_n)}(S(x_1,y_1,r_n))\\
& &+1_{[x_2\in S(x_1,y_1,r_n)]\cap[(x_1,x_2)\ not\ fails]}\ge
j_n\}\cap\{\mathcal{P}_{n(1-q_n)(1-v_n)}(S(x_2,y_2,r_n))\\
& &+1_{[x_1\in S(x_2,y_2,r_n)]\cap[(x_2,x_1)\ not\ fails]}\ge
j_n\}\big]\mathrm{d}y_1\mathrm{d}y_2\mathrm{d}x_2\mathrm{d}x_1\\
&\le&\frac{(1-v_n)^2n^2}{4\piup^2}\int_Q\int_{Q\cap
B(0,3r_n)}\int_0^{2\piup}\int_0^{2\piup}P\big[\{\mathcal{P}_{n(1-q_n)(1-v_n)}(S(0,y_1,r_n))\ge
j_n-1\}\\
& &\cap\{\mathcal{P}_{n(1-q_n)(1-v_n)}(S(x_2-x_1,y_2,r_n))\ge
j_n-1\}\big]\mathrm{d}y_1\mathrm{d}y_2\mathrm{d}x_2\mathrm{d}x_1\\
&\le&\frac{(1-v_n)^2n^2}{4\piup^2}\int_{B(0,3)}\int_0^{2\piup}\int_0^{2\piup}h_n(z,y_1,y_2)\mathrm{d}y_1\mathrm{d}y_2\mathrm{d}z
\end{eqnarray*}
where,
\begin{multline}
h_n(z,y_1,y_2):=P\big[\{\mathcal{H}_{nr_n^2(1-q_n)(1-v_n)}(S(0,y_1,1))\ge
j_n-1\}\\
\cap\{\mathcal{H}_{nr_n^2(1-q_n)(1-v_n)}(S(z,y_2,1))\ge
j_n-1\}\big].\nonumber
\end{multline}
By (\ref{3}), we choose $M\in\mathbb{N}$, such that
$j_n^M/\mu_n^{M+1}\rightarrow\infty$, as $n\rightarrow\infty$. Then
we have
$$P[\mathcal{H}_{nr_n^2(1-q_n)(1-v_n)}(S(0,y_1,1))\ge j_n+M]\le \xi_n\Big(\frac{\mu_n}{j_n}\Big)^M$$
and
$$P[\mathcal{H}_{nr_n^2(1-q_n)(1-v_n)}(S(0,y_1,1))\in\{j_n-1,j_n,\cdots,j_n+M-1\}]\le2\xi_n\frac{j_n}{\mu_n}$$
since
$P[\mathcal{H}_{nr_n^2(1-q_n)(1-v_n)}(S(0,y_1,1))=j_n-1]\le\xi_n(j_n/\mu_n)$
and $P[\mathcal{H}_{nr_n^2(1-q_n)(1-v_n)}(S(0,y_1\\,1))\ge
j_n]\le\xi_n(j_n/\mu_n)$ when $n$ is large enough.

Let $\eta_{z,y_1,y_2}:=2|S(0,y_1,1)\backslash S(z,y_2,1)|/\alpha$,
then we see that the conditional distribution of
$\mathcal{H}_{nr_n^2(1-q_n)(1-v_n)}(S(z,y_2,1))$, given that
$\mathcal{H}_{nr_n^2(1-q_n)(1-v_n)}(S(0,y_1,1))=j_n+M$ is the sum of
two independent random variables $j_n+M-U$ and $V$, where $U\sim
Bin(j_n+M,\eta_{z,y_1,y_2})$ and $V\sim Poi(\alpha
nr_n^2\eta_{z,y_1,y_2}(1-q_n)(1-v_n)/2)$. Provided $n$ is large
enough so that $M+1<j_n\eta_{z,y_1,y_2}/5$, if
$U>3j_n\eta_{z,y_1,y_2}/5$ and $V<j_n\eta_{z,y_1,y_2}/5$ then
$j_n+M-U+V<j_n-1$. Now, note that
$P[\mathcal{H}_{nr_n^2(1-q_n)(1-v_n)}(S(z,y_2,1))\ge
j_n-1|\mathcal{H}_{nr_n^2(1-q_n)(1-v_n)}(S(0,y_1,1))=k]$ is an
increasing function of $k$, and by Chernoff bounds, there exists a
constant $\beta>0$, for any $z\in B(0,3)$ and $n$ large enough, if
$\eta_{z,y_1,y_2}>5(M+1)/j_n$ then
\begin{eqnarray*}
\lefteqn{\max_{j_n-1\le k\le
j_n+M-1}P[\mathcal{H}_{nr_n^2(1-q_n)(1-v_n)}(S(z,y_2,1))\ge
j_n-1|\mathcal{H}_{nr_n^2(1-q_n)(1-v_n)}(S(0,y_1,1))=k]}\\
& &\le P[Bin(j_n+M,\eta_{z,y_1,y_2})\le
3j_n\eta_{z,y_1,y_2}/5]+P[Poi(\alpha
nr_n^2\eta_{z,y_1,y_2}(1-q_n)(1-v_n)/2)\ge j_n\eta_{z,y_1,y_2}/5]\\
& &\le2e^{-\beta j_n\eta_{z,y_1,y_2}},
\end{eqnarray*}
whereas if $\eta_{z,y_1,y_2}\le 5(M+1)/j_n$, then
$e^{5\beta(M+1)}e^{-\beta j_n\eta_{z,y_1,y_2}}\ge 1$. Take
$c_1=2\vee e^{5\beta(M+1)}$, for any $z\in B(0,3)$, we have
$$
\max_{j_n-1\le k\le
j_n+M-1}P[\mathcal{H}_{nr_n^2(1-q_n)(1-v_n)}(S(z,y_2,1))\ge
j_n-1|\mathcal{H}_{nr_n^2(1-q_n)(1-v_n)}(S(0,y_1,1))=k]$$
$$\le c_1e^{-\beta j_n\eta_{z,y_1,y_2}}.$$
Therefore, the above discussion gives
$$h_n(z,y_1,y_2)\le \xi_n\Big(\frac{\mu_n}{j_n}\Big)^M+2c_1\xi_n\frac{j_n}{\mu_n}e^{-\beta j_n\eta_{z,y_1,y_2}}$$
for $n$ large enough. Now since $\inf_{z\in
B(0,3)}\{\eta_{z,y_1,y_2}/||z||\}>0$ uniformly in $y_1$ and $y_2$,
there exists a constant $\gamma>0$ such that
\begin{eqnarray*}
\int_0^{2\piup}\int_0^{2\piup}\int_{B(0,3)}h_n(z,y_1,y_2)-
\xi_n\Big(\frac{\mu_n}{j_n}\Big)^M\mathrm{d}z\mathrm{d}y_1\mathrm{d}y_2
&\le&8\piup^2c_1\xi_n\big(\frac{j_n}{\mu_n}\big)\int_{B(0,3)}e^{-\gamma j_n||z||}\mathrm{d}z\\
&\le&c'\xi_n\big(\frac{j_n}{\mu_n}\big)\Gamma(2)/(\gamma j_n)^2.
\end{eqnarray*}
Accordingly, by the choice of $M$,
$$((1-v_n)n\xi_n)^{-1}I_{2,n}^{out}\le c((1-v_n)n\xi_n)^{-1}\cdot\mu_nn\xi_n\Big[\big(\frac{\mu_n}{j_n}\big)^M+\frac1{\mu_nj_n}\Big]\rightarrow0,$$
as $n\rightarrow\infty$. The out-degree case hereby follows from
(\ref{4}).

For in-degree, set
$$W_n^{'in}:=\sum_{i=1}^{N_{n(1-v_n)}}1_{[\#\{X_j\in\mathcal{P}_{n(1-v_n)(1-q_n)}|X_i\in S(X_j,Y_j,r_n)\}\ge j_n+1]\cap[X_i\in
Q]}$$ Then by Palm theory for Poisson process,
$$EW_n^{'in}\sim (1-v_n)n\int_QP[Poi\big(\frac{\alpha}{2}nr_n^2(1-v_n)(1-q_n)\big)\ge j_n]\mathrm{d}x=(1-v_n)n\xi_n,$$
as $n\rightarrow\infty$. By Lemma 1, we thereby obtain
\begin{equation}
\big|P(W_n^{'in}=0)-e^{-EW_n^{'in}}\big|\le
\min\Big(1,\frac{2}{(1-v_n)n\xi_n}\Big)(I_{1,n}^{in}+I_{2,n}^{in}),\label{5}
\end{equation}
for large enough $n$. We have
\begin{eqnarray*}
I_{1,n}^{in}&=&(1-v_n)^2n^2\int_{\mathbb{R}^2}P[\mathcal{\tilde{P}}_{n(1-q_n)(1-v_n)}(B(x_1,r_n))\ge
j_n]\\
&
&\int_{B(x_1,3r_n)}P[\mathcal{\tilde{P}}_{n(1-q_n)(1-v_n)}(B(x_2,r_n))\ge
j_n]f(x_2)\mathrm{d}x_2f(x_1)\mathrm{d}x_1\\
&\le&n^2\xi_n^2\piup(3r_n)^2
\end{eqnarray*}
Thus, by (\ref{3}), $((1-v_n)n\xi_n)^{-1}I_{1,n}^{in}\le
c\mu_n\xi_n\rightarrow0.$ On the other hand, we have,
\begin{eqnarray*}
I_{2,n}^{in}&=&(1-v_n)^2n^2\int_Q\int_{Q\cap
B(x_1,3r_n)}P\big[\{\mathcal{\tilde{P}}_{n(1-q_n)(1-v_n)}(B(x_1,r_n))\\
& &+1_{[x_1\in S(x_2,Y_2,r_n)]\cap[(x_2,x_1)\ not\ fails]}\ge
j_n\}\cap\{\mathcal{\tilde{P}}_{n(1-q_n)(1-v_n)}(B(x_2,r_n))\\
& &+1_{[x_2\in S(x_1,Y_1,r_n)]\cap[(x_1,x_2)\ not\ fails]}\ge
j_n\}\big]\mathrm{d}x_2\mathrm{d}x_1\\
&\le&(1-v_n)^2n^2\int_Q\int_{Q\cap
B(0,3r_n)}P\big[\{\mathcal{\tilde{P}}_{n(1-q_n)(1-v_n)}(B(0,r_n))\ge
j_n-1\}\\
& &\cap\{\mathcal{\tilde{P}}_{n(1-q_n)(1-v_n)}(B(x_2-x_1,r_n))\ge
j_n-1\}\big]\mathrm{d}x_2\mathrm{d}x_1\\
&\le&(1-v_n)^2n^2\int_{B(0,3)}g_n(z)\mathrm{d}z
\end{eqnarray*}
where,
\begin{multline}
g_n(z):=P\big[\{\mathcal{H}_{\frac{\alpha}{2\piup}nr_n^2(1-q_n)(1-v_n)}(B(0,1))\ge
j_n-1\}\\
\cap\{\mathcal{H}_{\frac{\alpha}{2\piup}nr_n^2(1-q_n)(1-v_n)}(B(z,1))\ge
j_n-1\}\big].\nonumber
\end{multline}
Also by (\ref{3}), we choose $M$ as above. Then we have
$$P[\mathcal{H}_{\frac{\alpha}{2\piup}nr_n^2(1-q_n)(1-v_n)}(B(0,1))\ge j_n+M]\le \xi_n\Big(\frac{\mu_n}{j_n}\Big)^M$$
and
$$P[\mathcal{H}_{\frac{\alpha}{2\piup}nr_n^2(1-q_n)(1-v_n)}(B(0,1))\in\{j_n-1,j_n,\cdots,j_n+M-1\}]\le2\xi_n\frac{j_n}{\mu_n}$$
since
$P[\mathcal{H}_{\frac{\alpha}{2\piup}nr_n^2(1-q_n)(1-v_n)}(B(0,1))=j_n-1]\le\xi_n(j_n/\mu_n)$
and
$P[\mathcal{H}_{\frac{\alpha}{2\piup}nr_n^2(1-q_n)(1-v_n)}(B(0,1))\\\ge
j_n]\le\xi_n(j_n/\mu_n)$ when $n$ is large enough. Let
$\delta_z:=|B(0,1)\backslash B(z,1)|/\piup$, hereby the conditional
distribution of
$\mathcal{H}_{\frac{\alpha}{2\piup}nr_n^2(1-q_n)(1-v_n)}(B(z,1))$,
given that
$\mathcal{H}_{\frac{\alpha}{2\piup}nr_n^2(1-q_n)(1-v_n)}(B(0,1))=j_n+M$
is the sum of two independent random variables $j_n+M-U$ and $V$,
where $U\sim Bin(j_n+M,\delta_z)$ and $V\sim Poi(\alpha
nr_n^2\delta_z(1-q_n)(1-v_n)/2)$. Provided $n$ is large enough so
that $M+1<j_n\delta_z/5$, if $U>3j_n\delta_z/5$ and
$V<j_n\delta_z/5$ then $j_n+M-U+V<j_n-1$. Now, note that
$P[\mathcal{H}_{\frac{\alpha}{2\piup}nr_n^2(1-q_n)(1-v_n)}(B(z,1))\ge
j_n-1|\mathcal{H}_{\frac{\alpha}{2\piup}nr_n^2(1-q_n)(1-v_n)}(B(0,1))=k]$
is an increasing function of $k$, and by Chernoff bounds, there
exists a constant $\beta>0$, for any $z\in B(0,3)$ and $n$ large
enough, if $\delta_z>5(M+1)/j_n$ then
\begin{eqnarray*}
\lefteqn{\max_{j_n-1\le k\le
j_n+M-1}P[\mathcal{H}_{\frac{\alpha}{2\piup}nr_n^2(1-q_n)(1-v_n)}(B(z,1))\ge
j_n-1|\mathcal{H}_{\frac{\alpha}{2\piup}nr_n^2(1-q_n)(1-v_n)}(B(0,1))=k]}\\
& &\le P[Bin(j_n+M,\delta_z)\le 3j_n\delta_z/5]+P[Poi(\alpha
nr_n^2\delta_z(1-q_n)(1-v_n)/2)\ge j_n\delta_z/5]\\
& &\le2e^{-\beta j_n\delta_z},
\end{eqnarray*}
whereas if $\delta_z\le 5(M+1)/j_n$, then $e^{5\beta(M+1)}e^{-\beta
j_n\delta_z}\ge 1$. Take $c_2=2\vee e^{5\beta(M+1)}$, for any $z\in
B(0,3)$, we have
$$
\max_{j_n-1\le k\le
j_n+M-1}P[\mathcal{H}_{\frac{\alpha}{2\piup}nr_n^2(1-q_n)(1-v_n)}(B(z,1))\ge
j_n-1|\mathcal{H}_{\frac{\alpha}{2\piup}nr_n^2(1-q_n)(1-v_n)}(B(0,1))=k]$$
$$\le c_2e^{-\beta j_n\delta_z}.$$
Consequently the above discussion gives
$$g_n(z)\le \xi_n\Big(\frac{\mu_n}{j_n}\Big)^M+2c_2\xi_n\frac{j_n}{\mu_n}e^{-\beta j_n\delta_z}$$
for $n$ large enough. Now since $\inf_{z\in
B(0,3)}\{\delta_z/||z||\}>0$, there exists a constant $\gamma>0$
such that
\begin{eqnarray*}
\int_{B(0,3)}g_n(z)- \xi_n\Big(\frac{\mu_n}{j_n}\Big)^M\mathrm{d}z
&\le&2c_2\xi_n\big(\frac{j_n}{\mu_n}\big)\int_{B(0,3)}e^{-\gamma j_n||z||}\mathrm{d}z\\
&\le&c'\xi_n\big(\frac{j_n}{\mu_n}\big)\Gamma(2)/(\gamma j_n)^2.
\end{eqnarray*}
Thus argue as the out-degree case,
$((1-v_n)n\xi_n)^{-1}I_{2,n}^{in}$ tends to $0$, as
$n\rightarrow\infty$. We hereby complete the proof by using
(\ref{5}). $\Box$
\medskip

Now we extend Proposition 1 from $\mathcal{P}_n$ to $\mathcal{X}_n$.

\medskip
\noindent\textbf{Proposition 2.}\itshape \quad Let
$\mu_n:=\frac{\alpha}{2}nr_n^2(1-v_n)(1-q_n)$. Suppose
$\inf_{n>0}\mu_n>0$ and
$\lim_{n\rightarrow\infty}\frac{\mu_n}{n^{1/6}}\\=0$. Suppose
$\{j_n\}_{n\ge 1}$ is an $\mathbb{N}-$valued sequence such that for
some $\varepsilon>0$, (\ref{3}) holds. Set $\xi_n:=P(Poi(\mu_n)\ge
j_n)$. Then
$$\lim_{n\rightarrow\infty}[P(\Delta_n^{out}<j_n)-e^{-n(1-v_n)\xi_n}]=0$$
and
$$\lim_{n\rightarrow\infty}[P(\Delta_n^{in}<j_n)-e^{-n(1-v_n)\xi_n}]=0.$$\normalfont

\medskip
\noindent\textbf{Proof}. For out-degree, we first assume that
$j_n\ge n^{1/5}$. We have
$$P(\Delta_n^{out}\ge n^{1/5})\le nP(Bin(n-1,\frac{\mu_n}{n})\ge n^{1/5})\rightarrow0,$$
as $n\rightarrow\infty$, by Chernoff bounds. Accordingly,
$P(\Delta_n^{out}\ge j_n)$ and $-n(1-v_n)\xi_n$ tend to zero. The
result then follows.

From now on, we thereby assume $j_n<n^{1/5}$ for all $n$, without
loss of generality. Set $\lambda_n:=n+n^{3/4}$, and let
$\mathcal{P}_{\lambda_n}$ be the Poisson point process coupled to
$\mathcal{X}_n$ with intensity $\lambda_nf$. Denote by
$\Delta_n^{+,out}$ the maximum out-degree in
$G_{\alpha}(\mathcal{P}_{\lambda_n},\mathcal{Y}_{N_{\lambda_n}},v_n,q_n,r_n)$.
Set
$\mu_n^+:=\frac{\piup\alpha}{2\piup}\lambda_nr_n^2(1-v_n)(1-q_n)$
and $\xi_n^+:=P(Poi(\mu_n^+)\ge j_n)$. Using Proposition 1 we have
\begin{equation}
\lim_{n\rightarrow\infty}[P(\Delta_n^{+,out}<j_n)-e^{-\lambda_n(1-v_n)\xi_n^+}]=0.\label{6}
\end{equation}
Since $1\le (\mu_n^+/\mu_n)^{j_n}=(1+n^{-1/4})^{j_n}\rightarrow1$
and $0\le\mu_n^+-\mu_n=n^{-1/4}\mu_n\rightarrow0$, as
$n\rightarrow\infty$, we get
$$1\le\frac{\xi_n^+}{\xi_n}\le\frac{e^{-\mu_n^+}(\mu_n^+)^{j_n}[1+\frac{\mu_n^+}{j_n}+(\frac{\mu_n^+}{j_n})^2+\cdots]}{e^{-\mu_n}\mu_n^{j_n}}\rightarrow1.$$
Then by setting $a_n=n(1-v_n)\xi_n$ and
$b_n=\frac{\lambda_n\xi_n^+}{n\xi_n}-1$, we have $b_n>0$ and
$b_n\rightarrow 0$. Observe that
$1-e^{-a_nb_n}\le1-e^{-\sqrt{b_n}}\rightarrow0$, if
$a_n\le1/\sqrt{b_n}$, while $e^{-a_n}<e^{-1/\sqrt{b_n}}\rightarrow0$
if $a_n>1/\sqrt{b_n}$. Consequently,
$$e^{-n(1-v_n)\xi_n}-e^{-\lambda_n(1-v_n)\xi_n^+}=e^{-a_n}(1-e^{-a_nb_n})\rightarrow 0,$$
as $n\rightarrow\infty$. Combining this with (\ref{6}), we have
$$\lim_{n\rightarrow\infty}[P(\Delta_n^{+,out}<j_n)-e^{-n(1-v_n)\xi_n}]=0.$$
Since
$P(\Delta_n^{+,out}<j_n)-P(\Delta_n^{out}<j_n)=P(\Delta_n^{+,out}<j_n\le\Delta_n^{out})-P(\Delta_n^{out}<j_n\le\Delta_n^{+,out}|n\le
N_{\lambda_n}\le n+2n^{3/4})\cdot P(n\le N_{\lambda_n}\le
n+2n^{3/4})-P(\Delta_n^{out}<j_n\le\Delta_n^{+,out}|\{N_{\lambda_n}<n\}\cup
\{N_{\lambda_n}>n+2n^{3/4}\})\cdot P(\{N_{\lambda_n}<n\}\cup
\{N_{\lambda_n}>n+2n^{3/4}\})$, and
$P(\Delta_n^{+,out}<j_n\le\Delta_n^{out})$ tends to 0, $P(n\le
N_{\lambda_n}\le n+2n^{3/4})$ tends to 1 by the Chebyshev
inequality, as $n\rightarrow\infty$, to prove the result it suffices
to prove that
$$\lim_{n\rightarrow\infty}P(\Delta_n^{out}<j_n\le\Delta_n^{+,out}|n\le
N_{\lambda_n}\le n+2n^{3/4})=0.$$ Now suppose $\Delta_n^{+,out}\ge
j_n$ and $n\le N_{\lambda_n}\le n+2n^{3/4}$, then there exists a
point in $\mathcal{P}_{\lambda_n}$ of out-degree at least $j_n$ in
$G_{\alpha}(\mathcal{P}_{\lambda_n},\mathcal{Y}_{N_{\lambda_n}},v_n,q_n,r_n)$.
Therefore
$$P[j_n>\Delta_n^{out}|\Delta_n^{+,out}\ge j_n,\ n\le N_{\lambda_n}\le n+2n^{3/4}]\le(j_n+1)\frac{2n^{3/4}}{n}\rightarrow0.$$
The out-degree case thereby follows by multiplication formula of
probability.

For in-degree, the same argument may be applied. $\Box$

\medskip
\noindent\textbf{Proof of Theorem 1}. Let $\xi_n(j):=P(Poi(\mu_n)\ge
j)$, then for $n\in\mathbb{N}$, take $j_n$ satisfying
$n\xi_n(j_n-1)>(1-v_n)^{-1}\ge n\xi_n(j_n)$. Set
$$
k_n:=\left\{\begin{array}{ll} j_n-1&,i\!f\
(1-v_n)n\xi_n(j_n)\le\sqrt{\frac{\xi_n(j_n)}{\xi_n(j_n-1)}}\\
j_n&,otherwise
\end{array}
\right.
$$
Take $\eta>0$ satisfying $(1+\varepsilon)^{-1}=1-2\eta$. Let
$i_n:=\lfloor\mu_n(\ln n)^{\eta}\rfloor$, then $i_n/(\ln
n)^{1-\eta}\rightarrow0$, as $n\rightarrow\infty$. Hence, by
Stirling formula,
$$
(1-v_n)n\xi_n(i_n)\ge(1-v_n)ne^{-1/12i_n}\frac{1}{\sqrt{2\piup
i_n}}e^{-i_n\ln(i_n/\mu_n)}\ge
cni_n^{-1/2}e^{-i_n\ln(i_n/\mu_n)}\rightarrow\infty.
$$
Thereby, $j_n\ge i_n$ and $j_n/\mu_n\rightarrow\infty$, as
$n\rightarrow\infty$. Hence, $\xi_n(j_n)/\xi_n(j_n-1)$,
$\xi_n(j_n+1)/\xi_n(j_n)$ and $\xi_n(j_n-1)/\xi_n(j_n-2)$ all tend
to zero. By the definition of $k_n$,
$(1-v_n)n\xi_n(k_n+1)\rightarrow0$ and
$(1-v_n)n\xi_n(k_n-1)\rightarrow\infty$, as $n\rightarrow\infty$.
Consequently, by Proposition 1, we have
$$P(\Delta_n^{'out}<k_n+1)\rightarrow1,\quad P(\Delta_n^{'out}<k_n-1)\rightarrow0,\quad P(\Delta_n^{'out}<k_n)-e^{-(1-v_n)n\xi_n(k_n)}\rightarrow0,$$
and
$$P(\Delta_n^{'in}<k_n+1)\rightarrow1,\quad P(\Delta_n^{'in}<k_n-1)\rightarrow0,\quad P(\Delta_n^{'in}<k_n)-e^{-(1-v_n)n\xi_n(k_n)}\rightarrow0.$$
Also by Proposition 2, we have
$$P(\Delta_n^{out}<k_n+1)\rightarrow1,\quad P(\Delta_n^{out}<k_n-1)\rightarrow0,\quad P(\Delta_n^{out}<k_n)-e^{-(1-v_n)n\xi_n(k_n)}\rightarrow0,$$
and
$$P(\Delta_n^{in}<k_n+1)\rightarrow1,\quad P(\Delta_n^{in}<k_n-1)\rightarrow0,\quad P(\Delta_n^{in}<k_n)-e^{-(1-v_n)n\xi_n(k_n)}\rightarrow0.$$
Thus, the proof is finally complete. $\Box$

\section{Open problems}

A natural question would be to ask what happens for other limiting
regime of $\mu_n$. We conjecture that when $\ln n\ll \mu_n\ll (\ln
n)^2$ and some regular conditions hold for fault probabilities,
there exist sequences $i_n$, $j_n$ such that for all $x$:
$$
P\Big(\frac{\Delta_n^{out/in}-i_n}{j_n}<x\Big)\rightarrow
e^{-e^{-x}}.
$$
Therefore, the focusing results will hold no longer in this case.
Results from extreme value theory suggest that might be the case. In
Erd\"os-R\'enyi random graph, a similar result holds \cite{15}.

Of course it would be of interest to consider the density function
other than the uniform one.

Note that our random faulty scaled sector graphs are still static
models, so what can be said about the behavior of maximum degrees of
a dynamic model? A direct and meaningful way to get a dynamic faulty
scaled sector graph is to give every point $X_i$ a random lifetime
$T_i$. Suppose that these lifetimes are independent random variables
with common distribution $F(t)=P(T_i\le t)$.



\end{document}